%



\input amstex
\magnification=\magstep1
\documentstyle{amsppt}
\UseAMSsymbols
\voffset=-3pc
\loadbold
\loadmsbm
\loadeufm
\UseAMSsymbols

\font\vastag=cmssbx10

\font\scVIII=cmcsc8

\font\ttVIII=cmtt8

\baselineskip=12pt
\parskip=6pt

\def\bC{\Bbb C}

\def\bR{\Bbb R}

\def\cE{\Cal E}

\def\cS{\Cal S}

\def\ep{\varepsilon}

\def\oU{\overline U}

\def\cln{{:}\;}
\def\id{\text{id}}

\def\dist{\text{ dist}}

\def\psc{\text{\rm{psc}}}

\def\ra{\rangle}
\def\la{\langle}
\NoBlackBoxes
\topmatter
\title PLURISUBHARMONIC DOMINATION IN BANACH SPACES
\endtitle
\author Imre Patyi\endauthor

\endtopmatter
\document

\subhead 1.\ Introduction\endsubhead

Given a complex manifold $M$, one says that plurisubharmonic, resp.\ holomorphic, domination is possible in $M$ if for any locally bounded function $u\cln M\to\bR$ there is a continuous plurisubharmonic function $w\cln M\to\bR$, resp.\ a Banach space $(V, \|\cdot\|_V)$ and a holomorphic function $h\cln M\to V$, such that
$$u(x)\le w(x),\quad\text{resp.}\quad u(x)\le\|h(x)\|_V, \quad \text{for every } x\in M.$$ 
These notions were introduced and studied by Lempert
in [L2]. The main result there is that if a Banach space $X$ has an unconditional basis and $\Omega \subset X$ is a pseudoconvex open set, then holomorphic, hence also plurisubharmonic domination is possible in $\Omega$. This result subsequently formed the basis for the study of analytic sheaves and cohomology groups in Banach spaces in [L3,LP,P1-2,S]. The goal of this paper is to prove that plurisubharmonic domination is possible in Banach spaces that have a Schauder basis (or are a direct summand of one that does, i.e., in separable Banach spaces with the bounded approximation property); this class includes all separable Banach spaces that occur in practice. In particular, domination is possible in the important Banach spaces $C[0,1]$ and $L^1[0,1]$, spaces that do not have an unconditional basis. More precisely, we shall prove 

\proclaim{Theorem 1.1}Let $X$ be a Banach space and $\Omega\subset X$ an open subset. If $X$ has a Schauder basis and $\Omega$ is pseudoconvex, then plurisubharmonic domination is possible in $\Omega$. The same holds if $X$ is just separable but $\Omega$ is convex.
\endproclaim

The second part of the theorem easily follows from the first. It would follow for all pseudoconvex $\Omega$ in a separable space if the following could be proved: Given a Banach space $X_0$, a closed subspace $X\subset X_0$, and a pseudoconvex $\Omega\subset X$, there is a pseudoconvex $\Omega_0\subset X_0$ such that $\Omega=\Omega_0\cap X$. --- It seems likely that under the assumptions of Theorem 1.1 holomorphic domination is also possible, but a proof of this will have to wait for another publication. When $\Omega$ is convex, holomorphic domination was already proved in [P3].

In its structure, the proof of Theorem 1.1 is similar to the proof in [L2]. The main new idea here and already in [P3] is that, at least when $\Omega=X$, it is possible to work through the proof with functions that are defined on all of $X$; and once the theorem is known for $\Omega=X$, the general case is not hard to prove. By contrast, [L2] treated all $\Omega\subset X$ on equal footing; it had to deal with holomorphic functions defined on subsets of $X$, and approximate them uniformly by entire functions. The required Runge--type approximation theorems are only known in Banach spaces with an unconditional basis (or, more generally, in spaces with a finite dimensional unconditional decomposition, see [J,L1,Me]) and this restricted the scope of [L2]. 

Since most results of [LP] depended on the hypothesis of plurisubharmonic domination, by Theorem 1.1 those results hold in spaces with a Schauder basis. For example, combining Theorem 1.1 with [LP, Theorem 2] gives the following generalization of Cartan's Theorems A and B:

\proclaim{Theorem 1.2}Suppose a Banach space $X$ has a Schauder basis, $\Omega\subset X$ is a pseudoconvex open subset, and $\cS\to\Omega$ a cohesive sheaf. Then:

{\rm(a)}\ There is a completely exact resolution $\dots\to\cE_1\to\cE_0\to\cS\to 0$; and

{\rm(b)}\ $H^q(\Omega,\cS)=0$, $q=1,2,\dots$, 
holds for the higher sheaf cohomology groups. 
\endproclaim

The notions occurring in this theorem are defined in [LP], to which we refer the reader; for fundamentals of complex analysis in Banach spaces the book [Mu] can be consulted.

\subhead 2.\ Ball bundles\endsubhead

Let $X$ be a Banach space. If $U\subset X$ is open, we write $\psc(U)$ for the family of continuous functions on $\oU$ that are plurisubharmonic on $U$. By saying that a subset $\Omega\subset X$ is pseudoconvex we imply that it is open. Suppose $X$ has a Schauder basis $e_1,e_2,\dots$, and introduce the projections $\pi_N\cln X\to X$, $$\pi_N\sum\limits^\infty_1\lambda_je_j=\sum\limits^N_1\lambda_je_j,\quad \lambda_j\in\bC;\qquad \pi_0=0,\quad \pi_\infty=\id.$$

We choose the norm $\|\cdot\|$ on $X$ so that for all $x\in X$
$$
\|\pi_Nx-\pi_Mx\|\le\|\pi_nx-\pi_mx\|,\quad 0\le n\le N\le M\le m\le\infty;\tag2.1
$$
thus $e_1,e_2,\dots$ form a bimonotone Schauder basis. Put furthermore $\varrho_N=\id-\pi_N$, $Y_N=\varrho_N X$. Given $N=0,1,2,\dots$, $A\subset\pi_NX\approx\bC^N$, and a continuous $r\cln A\to[0,\infty)$, the sets
$$\aligned
A(r)=&\{x\in X\cln  \pi_Nx\in A, \|\varrho_N x||<r(\pi_Nx)\}\qquad \text{and}\\
A[r]=&\{x\in X\cln  \pi_Nx\in A, \|\varrho_N x||\le r(\pi_Nx)\}\endaligned\tag2.2
$$
are ball bundles over finite dimensional bases. Any open $\Omega\subset X$ can be exhausted by such ball bundles as follows (see [L2, Section 3]). Let $d(x)=\min\{1,\dist(x,X\setminus\Omega)\}$ and, given $\alpha\in(0,1)$,
$$\aligned
D_N\la\alpha\ra&=\{t\in\pi_NX\cln \|t\|<\alpha N,\quad 1<\alpha N d(t)\},\\
\Omega_N\la\alpha\ra&=\{x\in X \cln \pi_N x\in D_N\la\alpha\ra, \quad \|\varrho_Nx\|<\alpha d (\pi_Nx)\}.\endaligned\tag2.3
$$

For example, if $\Omega=X$ then $\Omega_N\la\alpha\ra=\emptyset$ for $N\le 1/\alpha$ and
$$
\Omega_N\la\alpha\ra=\{x\in X\cln  \|\pi_N x\|<\alpha N, \|\varrho_N x\|<\alpha\}\qquad\text{for } N>1/\alpha.\tag2.4
$$

From now on we assume $\Omega$ is pseudoconvex.

\proclaim{Proposition 2.1}{\rm(a)}\ Each $\Omega_N\la\alpha\ra\subset\Omega$ is pseudoconvex.

{\rm(b)}\ ${\overline{\Omega_n\la\gamma\ra}}\subset\Omega_N\la\beta\ra$ if $n\le N, \gamma\le \beta/4$.

{\rm(c)}\ Given $\gamma$, each $x\in\Omega$ has a neighborhood contained in all but finitely many $\Omega_N\la\gamma\ra$.
\endproclaim

This is [L2, Proposition 3.1]. We also introduce another exhaustion of $\Omega$ by certain $\Omega^N\la\gamma\ra$; these are ball bundles with respect to the decomposition $X=\pi_{N+1}X\oplus Y_{N+1}$. Let $\gamma\in(0,1)$,
$$\aligned
p^N(s)=&\max\left\{{\|\pi_Ns\|\over N},{1\over N d(s)}, {\|\varrho_Ns\|\over d(s)}\right\}, \qquad s\in \Omega\cap\pi_{N+1}X,\\
D^N\la\gamma\ra=&\{s\in\Omega\cap\pi_{N+1}X\cln  p^N(s)<\gamma\},\\
\Omega^N\la\gamma\ra=&\{x\in\pi^{-1}_{N+1} D^N\la \gamma\ra \cln  \|\varrho_{N+1}x\|<\gamma d(\pi_{N+1}x)\}.
\endaligned \tag2.5$$

According to [L2, Proposition 3.2] we have:

\proclaim{Proposition 2.2}If $\gamma<1/4$ then $\Omega_N\la\gamma\ra\subset\Omega^N\la 4\gamma\ra$ and $\Omega^N\la\gamma\ra\subset\Omega_N\la 4\gamma\ra$.
\endproclaim

We shall also need the following analogs of [L2, Lemma 4.1, Proposition 4.2]:

\proclaim{Lemma 2.3}Suppose $A_2\subset\subset A_3\subset\subset A_4$ are relatively open subsets of $\pi_NX, A_1\subset A_2$ is compact and holomophically convex in $A_4$. Let $r_i\cln A_4\to (0,\infty)$ be continuous, $i=1,2,3$, $r_i< r_{i+1}$, and $-\log r_1$ plurisubharmonic. Given $v\in\psc(X)$, there is a continuous plurisubharmonic $w\cln \pi_N^{-1}A_4\to \bR$ such that
$$
w(x)\cases <0, &\text{if\quad $x\in A_1[r_1]$}\\
           >v(x), &\text{if\quad $x\in A_3(r_3)\setminus A_2(r_2)$}.\endcases
	    $$
\endproclaim

\demo{Proof}As in the proof of [L2, Lemma 4.1] we construct a Banach space $(V, \|\cdot\|_V)$ and a holomorphic function $\psi\cln \pi^{-1}_NA_4\to V$ such that $\|\psi\|_V<1/4$ on $A_1[r_1]$ and $\|\psi\|_V>4$ on $A_3(r_3)\setminus A_2(r_2)$. (This corresponds to choosing $r_4=\infty$ there. The construction does not use the approximation hypothesis of Lemma 4.1.) Since $v$ is bounded on a neighborhood of the compact set $A_1$, there is a $q\in(0,\infty)$ such that
$$
v(y)\le q,\qquad \text{if\quad $\pi_Ny\in A_1$}\quad \text{and $\|\varrho_Ny\|\le 4^{-q}\,\max\limits_{A_1}\, r_1$}.\tag2.6
$$

Let $K$ be the set of linear forms on $V$ of norm $\le 1$, and define the continuous plurisubharmonic function $w$ by
$$
w(x)=2q\log(\|\psi(x)\|+1/4)+\sup\limits_{k\in K} v(\pi_N x+(k\psi(x))^q\varrho_Nx), \quad x\in\pi^{-1}_N A_4.\tag2.7
$$

	To check that $w$ is continuous and plurisubharmonic, it is enough to do
	the same for the sup in (2.7).
	As the argument $a(x,k)=\pi_Nx+(k\psi(x))^q\varrho_Nx$ is a holomorphic function
	of $(x,k)\in\pi^{-1}_N(A_4)\times V^*$ and $v$ is continuous and plurisubharmonic
	on $X$, it is enough to show that the supremum is is finite and continuous.
	To that end it is enough to show that if $x$ is confined to a compact set
	$C\subset\pi^{-1}_N(A_4)$ and $k$ to a bounded set (such as $K$), then
	the argument $a(x,k)$ is also confined to a compact set, and that $a(x,k)$ is
	Lipschitz continuous on any set $(x,k)\in B\times K$, where 
	$B\subset\pi^{-1}_N(A_4)$ is any open set on which $\psi$ is Lipschitz continuous
	(e.g., any ball small enough about any point $x$).
	Both these claims are obvious.

If $x\in A_1[r_1]$, then $|k\psi(x)|<1/4$ and $\|\psi(x)\|+1/4<e^{-1/2}$, hence $w(x)<-q+q=0$ by (2.6). If $x\in A_3(r_3)\setminus A_2(r_2)$ then the first term in (2.7) is positive, and there is a $k\in K$ such that $k\psi(x)=1$; therefore $w(x)>v(\pi_Nx+\varrho_Nx)=v(x)$, q.e.d. 
\enddemo

\proclaim{Proposition 2.4}Let $0<4^2\beta<\alpha<4^{-2}, N=1,2,\dots$. If $v\in\psc(X)$, there is a continuous plurisubharmonic $w\cln \pi^{-1}_{N+1}\Omega\to \bR$ such that
$$
w(x)\cases <0, &\text{if\quad $x\in \Omega_N\la\beta\ra$}\\
            >v(x), &\text{if\quad $x\in \Omega_{N+1}\la 4\alpha\ra\setminus\Omega_N \la \alpha\ra$}.\endcases
$$
\endproclaim

\demo{Proof}Let $A_4=\Omega\cap\pi_{N+1}X$ and with notation in (2.3), (2.5) define bounded sets 
$$A_1=\{s\in A_4\cln  p^N(s)\le 4\beta\},\quad A_2=D^N\la\alpha/4\ra,\quad A_3=D_{N+1}\la4\alpha\ra,$$
of which $A_1$ is compact, and $A_2, A_3$ are open in $A_4$. Let furthermore $r_1=4\beta d$, $r_2=\alpha d/4$, and $r_3=4\alpha d$. We apply Lemma 2.3 with $N$ replaced by $N+1$. Clearly $A_1\subset A_2$ is plurisubharmonically, hence holomorphically convex in $A_4$ (see [H, Theorem 4.3.4]). By (2.3), (2.5), and by Proposition 2.2
$$
A_1[r_1]\supset\Omega^N\la 4\beta\ra\supset \Omega_N\la\beta\ra, \qquad A_2[r_2]=\Omega^N\la \alpha/4\ra\subset \Omega_N\la\alpha\ra.
$$
Proposition 2.1b implies $\overline{A_2(r_2)}\subset\Omega_{N+1}\la 4\alpha\ra=A_3(r_3)$. Intersecting with $\pi_{N+1}X$, $\overline{A_2}\subset A_3$ follows, and $\overline{A_3}\subset A_4$ is obvious. Therefore by Lemma 2.3 there is a continuous plurisubharmonic $w\cln \pi^{-1}_{N+1} A_4=\pi^{-1}_{N+1}\Omega\to \bR$ as claimed: $w<0$ on $A_1[r_1]\supset\Omega_N\la\beta\ra$ and $w>v$ on $A_3(r_3)\setminus A_2(r_2)\supset \Omega_{N+1}\la 4\alpha\ra\setminus\Omega_N\la\alpha\ra$.
\enddemo

\subhead 3.\ Domination in the whole space\endsubhead

We prove the following special case of Theorem 1.1:

\proclaim{Proposition 3.1}Suppose a Banach space $X$ has a Schauder basis and $u\cln  X\to\bR$ is a locally bounded function. There is a $w\in\psc(X)$ such that $u(x)<w(x)$ for $x\in X$.
\endproclaim

We shall use the assumptions and the notation of section 2. If $x\in X$ and $\ep>0$, $B(x,\ep)\subset X$ will stand for the ball of radius $\ep$, centered at $x$. The key is the following

\proclaim{Proposition 3.2}Given $u\cln X\to\bR$, suppose there is an $\ep>0$ and for every $x\in X$ a $w_x\in\psc(X)$ such that $u<w_x$ on $B(x,\ep)$. Then there is a $w\in\psc(X)$ such that $u<w$.
\endproclaim

\demo{Proof}We can assume $u>0$ everywhere. Fix a positive $\alpha<\min(\ep,4^{-2})$, and with $N=1,2,\dots$, consider the compact set $A=\overline{\Omega_N\la\alpha\ra}\cap\pi_NX$; here $\Omega_N\la\alpha\ra$ refers to the exhaustion of $X=\Omega$ defined in (2.3) or (2.4). As each $t\in A$ has a neighborhood $U\subset\pi_NX$ such that $\Omega_N\la\alpha\ra\cap\pi^{-1}_NU\subset B(t,\ep)$, there is a finite $T\subset A$ such that
$$\Omega_N\la\alpha\ra\subset\bigcup_{t\in T}B(t,\ep).$$
It follows that $v_N=\max\{w_t\cln t\in T\}\in\psc(X)$ satisfies $v_N>u$ on $\Omega_N\la\alpha\ra$. Let $0<\beta<\alpha/4^2$. By Proposition 2.4, there is $w_N\in\psc(X)$ with
$$
w_N(x)\cases <0, &\text{if\quad $x\in \Omega_N\la\beta\ra$}\\
            >v_N(x), &\text{if\quad $x\in \Omega_{N+1}\la\alpha\ra\setminus \Omega_N\la\alpha\ra$},\endcases
$$
and $w=\sup\{v_1,w_1, w_2,\dots\}$ will do.
\enddemo

\demo{Proof of Proposition 3.1}Suppose the claim is not true, and $u$ cannot be dominated by any $w\in\psc(X)$. In light of Proposition 3.2 there must be a ball $B(x_1,1)$ on which $u$ cannot be dominated by a $w\in\psc(X)$, i.e.,
$$
u_1=\cases u &\text{on $B(x_1,1)$}\\
            0 &\text{on $X\setminus B(x_1,1)$},\endcases
$$
cannot be plurisubharmonically dominated. Again by Proposition 3.2, there must be a ball $B(x_2, 1/2)$ on which $u_1$ cannot be dominated, and so on. We obtain a sequence $B(x_k, 1/k)$ of balls such that
$$
u_k=\cases u &\text{on $\bigcap^k_1B(x_j,1/j)$}\\
            0 &\text{on $X\setminus \bigcap^k_1B(x_j,1/j)$},\endcases
$$ 
cannot be plurisubharmonically dominated. In particular $\bigcap^k_1B(x_j, 1/j)\ne\emptyset$. Hence $\|x_j-x_k\|<(1/j)+(1/k)$, and the $x_j$ have a limit $x$. But $u$ is bounded on some neighborhood of $x$, so on some $B(x_k, 1/k)$; hence $u_k$ is bounded and can be dominated by a constant. This is a contradiction, which then proves the claim.
\enddemo

\subhead 4.\ Domination in a general $\Omega$\endsubhead

Consider a pseudoconvex subset $\Omega$ of a Banach space $X$ that has a Schauder basis.

\proclaim{Proposition 4.1}Given a locally bounded $u\cln \Omega\to\bR$, there is a continuous plurisubharmonic $w\cln \Omega\to\bR$ such that $u(x)<w(x)$ for $x\in\Omega$.
\endproclaim

\demo{Proof}Again we make the assumptions and use notation introduced in Section 2. Fix $0<\alpha<4^{-2}$ and $0<\beta<\alpha/4^2$. For $N=0,1,\dots$ let $U_N=\bigcap_{j\ge N}\Omega_j\la 4\alpha\ra$. By Proposition 2.1c, these are open sets and exhaust all of $\Omega$. We prove by induction that there are $w_N\in\psc (U_{N+1})$ such that
$$
w_N\cases <0 &\text{on $\Omega_N\la\beta\ra$}\\
           >u &U_{N+1}\setminus U_N,\endcases
\qquad \text{and }\quad w_N>w_{N-1}\quad\text{on}\quad \partial U_N.\tag4.1
$$ 
(When $N=0$, the last requirement is vacuous, $U_0=\Omega_0\la 4\alpha\ra=\emptyset$.) The functions
$$
u_N=\cases u &\text{on $\overline{U_{N+1}}$}\\
            0 &\text{on $X\setminus \overline{U_{N+1}}$},\endcases
$$
are locally bounded. Applying Proposition 3.1 we obtain $w_0\in\psc(U_1)$ with $w_0>u_0$; then (4.1) is satisfied for $N=0$.

Next suppose that $w_0,\dots, w_{N-1}$ have already been found. Again by Proposition 3.1 there is $v\in\psc (X)$ such that $v>u_N$ on $X$ and $v>w_{N-1}$ on $\partial U_N$. Further, by Proposition 2.4 there is a continuous plurisubharmonic $v'\cln \pi_{N+1}^{-1}\Omega\to\bR$ such that
$$
v'\cases <0 &\text{on $\Omega_N\la\beta\ra$}\\
           >v &\text{on $\Omega_{N+1}\la 4\alpha\ra\setminus \Omega_N\la \alpha\ra$}.\endcases
$$ 
In view of Proposition 2.1b, $U_N\supset\Omega_N\la\alpha\ra$, and so $U_{N+1}\setminus U_N\subset\Omega_{N+1}\la 4\alpha\ra\setminus\Omega_N\la \alpha\ra$. It follows that $w_N=v'|U_{N+1}\in\psc (U_{N+1})$ satisfies (4.1).

Define $w\cln \Omega\to\bR$ by
$$
w(x)=\sup\{w_N(x), w_{N+1}(x),\dots\},\quad\text{if}\quad x\in U_{N+1}\setminus  U_N,\; N=0,1,2,\dots.\tag 4.2$$
By (4.1) and Proposition 2.1c, the sup is locally finite, and so defines a continuous plurisubharmonic function on $U_{N+1}\setminus\overline{U_N}$, hence on $\Omega\setminus\bigcup_{N\ge 1}\partial U_N$. But $w$ is also continuous and plurisubharmonic in some neighborhood of any $x_0\in\partial U_N$. Indeed, choose $N\le M$ so that
$$x_0\in\partial U_N\cap\partial U_{N+1}\cap\dots\cap\partial U_M\cap U_{M+1},\; \text{and}\quad x_0\notin\overline{U_{N-1}}.$$
By (4.1), $w_M(x_0)>w_{M-1}(x_0)>\dots>w_{N-1}(x_0)$. By continuity, it follows that for $x$ near $x_0$,
$$w(x)=\sup\{w_M(x), w_{M+1}(x),\dots\},\qquad\text{cf}\; (4.2).$$
Since $w_j$ for $j\ge M$ is continuous and plurisubharmonic on a neighborhood of $x_0$, so is $w$. Finally, (4.1) implies $w>u$, and the proof is complete.
\enddemo

\subhead 5.\ Separable spaces\endsubhead

Proposition 4.1 represents the first part of Theorem 1.1. To prove the second part, let $X$ be separable and $\Omega\subset X$ convex and open. 
Embed $X$ linearly into the space $X_0=C[0,1]$, so that $X\subset X_0$ is a (closed) linear subspace. 
We can assume $0\in\Omega$. Let $B\subset X_0$ be an open ball centered at 0 such that $B\cap X\subset\Omega$. The convex hull $\Omega_0$ of $B\cup\Omega$ is a convex, open subset of $X_0$. 
We claim that $\Omega_0\cap X=\Omega$. Indeed, suppose $p\in X\setminus\Omega$. By the Hahn--Banach separation theorem, there is a real linear form $f\cln X\to\bR$ such that
$$f(p)>f(x)\qquad\text{for}\quad x\in\Omega.\tag 5.1$$
In particular, $f(p)>f(x)$ for $x\in B\cap X$. If $f_0\cln X_0\to\bR$ denotes a linear extension of $f$, having the same norm as $f$, then
$$f_0(p)>f_0(x)\qquad\text{for}\quad x\in B.\tag 5.2$$
(5.1) and (5.2) imply $f_0(p)>f_0(x)$ for $x\in\Omega_0$, whence $p\notin\Omega_0$ as claimed. It follows that $\Omega=\Omega_0\cap X$ is closed in $\Omega_0$.

Any locally bounded $u\cln \Omega\to\bR$ extends to the locally bounded function
$$
u_0=\cases u &\text{on $\Omega$}\\
           0 &\text{on $\Omega_0\setminus \Omega$}.\endcases
	    $$
Since $X_0=C[0,1]$ has a Schauder basis, by Proposition 4.1 $u_0$ can be dominated by a continuous plurisubharmonic $v_0$; then $v=v_0|\Omega$ dominates $u$, q.e.d.

\vskip0.20truein
\centerline{\vastag*~***~*}
\vskip0.15truein
{\scVIII  
        Imre Patyi,
        Department of Mathematics and Statistics,
        Georgia State University,
        Atlanta, GA 30303-3083, USA,
        {\ttVIII ipatyi\@gsu.edu}
}
\enddocument
\bye